\documentclass[12pt]{article}

\usepackage{amsmath, amsthm, amsopn, amssymb, enumerate, fullpage}

\usepackage{geometry}                
\usepackage{graphicx}
\usepackage{amssymb}
\usepackage{epstopdf}
\newtheorem{theorem}{Theorem}[section]

\newtheorem{problem}[theorem]{Problem}
\newtheorem{proposition}[theorem]{Proposition}

\newtheorem{claim}[theorem]{Claim}

\title{Partitioning the hypercube into smaller hypercubes}
\author{Noga Alon, J\'ozsef Balogh, VP}

\author{%
Noga Alon\footnote{Princeton University, Princeton, NJ 08544,
USA and Tel Aviv University, Tel Aviv 69978, Israel. E-mail:
\texttt{nalon@math.princeton.edu}. Research supported in part by NSF grant
DMS-2154082.}
\and 
J\'ozsef Balogh\footnote{Department of Mathematics, 
University of Illinois at Urbana-Champaign, Urbana, Illinois 61801, 
USA. E-mail: \texttt{jobal@illinois.edu}. 
Research is partially supported by NSF Grant DMS-1764123,  
Arnold O. Beckman Research Award (UIUC Campus Research Board RB 24012), 
and the Langan Scholar Fund (UIUC).}
\and
Vladimir N. Potapov\footnote{E-mail: \texttt{quasigroup349@gmail.com.}}
}

\begin{document}
\maketitle

\begin{abstract}
Denote by $Q_d$ the $d$-dimensional hypercube. 
We estimate the number of ways the vertex set of $Q_d$ 
can be partitioned into 
vertex disjoint smaller cubes. Among other results, we prove
that the asymptotic order
of this function is larger than the number of 
perfect matchings of $Q_d$ by an exponential factor in the number
of vertices, and not by a larger factor. We also  describe 
and address several new (and old) related questions.
\end{abstract}

\section{The problem and main results}

Denote by $Q_d:=\{0,1\}^d$ the $d$-dimensional hypercube, and let $f(d)$
be  the number  of partitions (tilings) of the vertex set of $Q_d$
in which each of the classes spans a  subhypercube. Let
$f_{S}(d)$ denote the number of such partitions where each of the classes
has dimension in $S\subset\{0,1,\ldots,d\}$. In particular, $f_{\le 2}(d)$
counts the partitions where the parts are singletons, edges or spanning
a $2$-cube.\\ For a graph $G$, denote by $m(G)$ the number of perfect matchings in $G$, and by $m'(G)$ the number of matchings, and write  $m(d):=m(Q_d)=f_1(d)$ and $m'(d)=m'(Q_d)=f_{\le 1}(d)$.

\smallskip

An easy observation is that the number of (perfect) matchings of
$Q_d$ is a lower bound on $f(d)$. Determining or estimating the
number of  (perfect) matchings of graphs is a classical problem,
see, e.g., \cite{AF,csikvari,friedland,JP}. 
Although it has been studied for $Q_d$,
besides determining it for small $d$, only general results are known about
its asymptotics, as we discuss later. Below we collect many
problems that might be of interest; and  we will address several of them.

\begin{problem}\label{all}
Determine or estimate the functions \\
(i) $m(d)$, \quad\quad \quad\quad (ii)  $m'(d)$, 
\quad \quad \quad \quad  (iii) $f(d)$,  
\quad\quad \quad\quad (iv)  $f_{\le 2}(d).$
\end{problem}

Problem (i) is natural and was studied in the literature, (ii) is
a natural extension. Problem (iii) was raised by Gadouleau at
the 29th British Combinatorial Conference in  2022, motivated by the paper of 
 Bridoux,  Gadouleau, and
Theyssier~\cite{BGT20}. Problem~(iv)
is a variant we suggest here, and seems 
to be closely related to (iii).

As pointed out in \cite{BGT20},
the first few terms of $f(d)$ are given in the
Online Encyc\-lopedia of Integer Sequences (OEIS) 
as A018926. Starting from $d=0$ these are:
$1, 2, 8, 154, 89512, $ $ 71319425714$.

Graham and Harary~\cite{harary} studied the number of 
perfect matchings of the hypercube. They
calculated\\
$
\quad\quad \ \ \hskip 1cm m(1)=1,\quad m(2)=2, \quad m(3)=9, \quad m(4)=272, \quad m(5)=589185.
$

Using dynamic programming, \"Osterg\aa rd and Pettersson~\cite{oster} 
determined $m(6)$ (which was also known earlier) and $m(7)$:\\
$m(6)=16 332 454 526 976, \quad      \quad     
m(7) = 391 689 748 492 473 664 721 077 609 089.
$

Determining precisely, or even providing satisfactory estimates for 
the functions $m(d),$ $ m'(d)$ and $f_S$ for some interesting families $S$
seems to be difficult, hence it would also be interesting to compare
the orders of magnitude of some of these functions.  
 Set
 \vglue -0.6cm
$$n:=2^{d-1}\quad\quad  \text{and}  \quad\quad N:=(d/e)^n.$$
Note that $n$ is the number of vertices in each vertex class of 
$Q_d$, and as mentioned below, $N$ is a rough estimate for the
number of perfect matchings in it.
\vspace{0.2cm}

\noindent 
The following
hierarchy follows from the definitions:
\vglue -0.4cm

\begin{equation}\label{chain}
m(d)   \    \le \                      m'(d)            \        \le    \        f_{\le 2}(d) \            \le \           f(d).
\end{equation}

We prove that the ratio $f(d)/m(d)$ is  exponential in $n$, and
wonder in which of the inequalities in \eqref{chain}
there is an exponential ratio.

\begin{problem}\label{ratio} What is the order of magnitude of the 
following ratios? In particular, which one  of them is an 
exponential function of $n$? \\
(i)   $m'(d)/m(d)$, \quad\quad (ii) $f_{\le 2}(d)/m'(d)$,  
\quad\quad (iii)  $f(d)/f_{\le 2}(d),$ \quad\quad \ (iv) $f(d)/m(d)$.
\end{problem}

Our  results and conjectures are summarized below.

Classical results easily imply Proposition~\ref{perfectmatching}, which
determines the main term of the asymptotics of $m(d)$, partially
solving   Problem~\ref{all}~(i).\footnote{The proof of
the upper bound, and that of a slightly weaker lower bound on $m(d)$ 
as stated in
Proposition~\ref{perfectmatching}, are posted in several class notes at
various course websites, see also \cite{LP}.
Here we provide a short proof, for completeness.}
Proposition~\ref{matching}  partially  solves   Problem~\ref{all}
(ii) and Problem~\ref{ratio} (i), showing that $m'(d)/m(d)$
is subexponential.
For both statements we state general results about regular 
bipartite graphs, which imply the
required bounds for the hypercube.
Proposition~\ref{expgap} shows that $f_{\le 2}(d)/m'(d)$ 
and $f_{\leq 2}(d)/m(d)$ are (at least) exponential
in $n$,  addressing Problem~\ref{ratio} (ii) and (iv).
In Proposition~\ref{main} we show that the ratio $f(d)/m(d)$ is at most
exponential in $n$, and a modification of its proof implies
that allowing to use
$2$-dimensional subcubes in the
partitions has a large impact on their numbers, see
Proposition~\ref{approx}.

Our final result, Proposition~\ref{ired},  shows that there are many {\it
irreducible tight} partitions. The precise definition of this notion appears 
before the statement of the result. 

We proceed with the formal statements of all the results.

\begin{proposition}\label{perfectmatching}
(i) Let $G$ be an $a$-regular bipartite graph with vertex classes of 
size $b$, where $a \geq 2$.  Then
\begin{equation}\label{prp13i}
 e^{b/2a}\cdot \left(\frac{a}{e}\right)^b \le m(G)\le 
  [\sqrt{2\pi a}\cdot 
e^{1/(12a)}]^{b/a}  \cdot \left(\frac{a}{e}\right)^b =
 2^{(1+o(1))\cdot \log a \cdot (b/2a)}\cdot \left(\frac{a}{e}\right)^b,
 \end{equation}
 when both $a,b\to \infty.$
 
 (ii) In particular for $Q_d$ we have 
 \begin{equation}\label{prp13ii}
 e^{n/2d}\cdot N\le m(d)\le 2^{(1+o(1))\cdot \log d\cdot (n/2d)}\cdot N.
 \end{equation}
\end{proposition}

The following result is addressing Problem~\ref{all} (ii) and
Problem~\ref{ratio} (i).
The lower bound in Proposition~\ref{matching} (i) also follows from 
known results in \cite{csikvari}, and the upper bound from the ones
in \cite{JP}. In particular, these results imply that 
$m'(G)/m(G)=e^{(2+o(1))n/\sqrt d}.$
Here we describe a short and
simple proof,  which does not provide
the optimal  constant in the exponent
obtained in the upper bound.

\begin{proposition}\label{matching}
(i) Let $G$ be an $a$-regular bipartite graph with vertex classes of
size $b$, where $a\ge 1$ and $b\to \infty.$ Then
\begin{equation}\label{prp14}  
m'(G)=m(G) \cdot 2^{\Theta(b/\sqrt a)}=\left(\frac{a}{e}\right)^b \cdot 2^{\Theta(b/\sqrt a)}.
\end{equation}
(ii) In particular for $Q_d$ we have 
\begin{equation}\label{prp144}  
m'(d)=m(d) \cdot 2^{\Theta(n/\sqrt d)}=N \cdot 2^{\Theta(n/\sqrt d)}.
\end{equation}
\end{proposition}

Next, we show that the number of cube partitions 
of $Q_d$ is exponentially larger than the number of matchings.

\begin{proposition}\label{expgap}
There exists a constant $c>1$ such that for all $d \geq 3$
\begin{equation}\label{prp15} 
c^{n}\cdot m'(d)\ \leq \ f_{\le 2} (d) \ \leq \ f(d).
\end{equation}
\end{proposition}

The following result is a simple upper bound on $f(d)$. Note 
that a sketch of the proof appeared already in~\cite{Po}.

\begin{proposition}\label{main}
$$  f(d) \le 
(d+1)^n ~\leq e^{n+n/d} \cdot N.
$$
\end{proposition}

The following Propositions show the effect of allowing $1$- and
$2$-dimensional (and in general small dimensional) subcubes in a partition.

\begin{proposition}\label{new}
For every fixed $r \geq 2$
$$   
N^{r/2^{r-1}-o(1)} \leq f_{0,r}(d) \leq 
f_{0, r,r+1,r+2,\ldots} (d) \le 
N^{r/2^{r-1}+o(1)}.$$
\end{proposition}

\begin{proposition}\label{approx}
$$   f_{0,1,  3,4,\ldots} (d) \le  \exp(20 n/d^{1/4})\cdot 
N \le \exp(20 n/d^{1/4}) \cdot f_{1}(d)=\exp(20 n/d^{1/4}) \cdot m(d).$$
\end{proposition}

It is natural to ask what happens if only $2$-dimensional
subcubes are  allowed to be in a partition. We believe that the same
bound as in Proposition~\ref{new}, with $r=2$ holds.
\begin{problem}\label{2dim}
Determine the asymptotic behaviour of the number $   f_{2} (d).$
\end{problem}

Additionally, we left open another question addressing 
Problem~\ref{ratio} (iii).

\begin{problem}\label{app012} Is it true that 
$   f (d)/ f_{\le 2}(d)$ is subexponential in $n$?
\end{problem}

A subcube partition is an {\it irreducible}, 
if  no subcube is spanned by a subfamily of
the partition, and it is  {\it tight} if  the partition is `proper',
in the sense that
every coordinate is used, i.e.,  every coordinate is fixed
in at least one subcube of the partition.
For example, the partition of the $2$-cube into the two subcubes
$\{0*\}$, in which  the first coordinate is fixed to $0$,
and $\{1*\}$, in which the first coordinate is fixed to $1$, is not a
tight partition, as the second coordinate is not used.

Irreducible subcube partitions appear in a work of Kullmann and
Zhao~\cite{KZ16}
and variants are described in several other
papers, see \cite{irreducible} for relevant references.
Note that it is not immediately clear that there exists an 
irreducible tight partition for every given (large) $d$.

Peitl and Szeider~\cite{PS23} enumerated all tight irreducible subcube
partitions for $ d = 3, 4$, and asked whether there are infinitely
many such partitions. Filmus, Hirsch,  Kurz,
Ihringer,  Riazanov,  Smal, and  Vinyals~\cite{irreducible} answered
this question in the affirmative, giving many explicit constructions of
tight irreducible subcube partitions. Here we prove the existence of many
more tight irreducible partitions.

\begin{proposition}\label{ired}
The number of irreducible tight 
partitions of the $d$-dimensional cube is at least
$c^{n}$ for some absolute constant $c>1$, where, as before, $n=2^{d-1}$.
\end{proposition}

\noindent We shall frequently use the Stirling formula:
$$ \sqrt{2\pi m}\cdot  \left(\frac{m}{e}\right)^m\ \le\  m! \ \le \  e^{1/(12m)}\cdot \sqrt{2\pi m}\cdot \left(\frac{m}{e}\right)^m.$$
Furthermore, we use the binary entropy function 
$h(x):=-x\log_2x -(1-x)\log_2(1-x)$ to estimate the 
binomial coefficients for $0<x<1/2$:

\begin{equation}\label{entropy}
\sum_{k\le xn} \binom{n}{k} \ =\  \Theta(1) \cdot  \binom{n}{xn} =2^{(1+o(1))h(x)\cdot n}.
\end{equation}

\noindent All logarithms throughout the paper 
have base $2$, unless otherwise indicated. To simplify the
presentation, we omit all floor and ceiling signs whenever these
are not crucial. We also assume, whenever this is needed, that
$d$ (and hence also $n$) is sufficiently large.

In Section~\ref{pm} we prove Propositions~\ref{perfectmatching},
\ref{matching}, and \ref{expgap}.  In Section~\ref{uppfd}
we prove  Propositions~\ref{main},~\ref{new},~\ref{approx}, and~\ref{ired}.

\vskip 0.1 cm

\noindent
{\bf Remark:} We mention very briefly two  
motivations for the study of $f(d)$.
First, $f(d)$ is the number of so-called bijective commutative Boolean
networks of dimension $d$.  See \cite{BGT20}  for the definition of this
notion and the fact that there is a bijection 
between partitions of the $d$-cube into subcubes and such networks.

Secondly, $f(d)$ is the number of instances of SAT on $d$ Boolean
variables such that any truth assignment fails to satisfy exactly one
clause of the instance. Indeed, to any clause (e.g. $x_1 \lor \neg x_3$)
we can associate the subcube of assignments that fail to satisfy it
($x_1 = 0, x_3 = 1$ in the example above). Then those subcubes partition the
$d$-cube if and only if any truth assignment belongs to exactly one of
them, i.e.,~it fails to satisfy exactly one clause.  This is equivalent
to the description from OEIS, which reads: 
``the number of  ways to make a tautology from disjoint terms with $d$
Boolean variables''.

\section{Matchings and $2$-dimensional subcubes}\label{pm}


\noindent

First, we list some of the classical results 
that we shall use in our proofs.

{\bf Bregman-Minc inequality:} \ The celebrated Bregman-Minc inequality, conjectured by Minc~\cite{minc}
and proved by Bregman~\cite{bregman} (see also \cite{Sch}, \cite{as},
\cite{Rad} for short proofs)
implies that the maximum 
possible number of
perfect matchings in an $a$-regular bipartite graph with $b$
vertices in each vertex class is at most $(a!)^{b/a}$. Equality is
achieved when $b$ is divisible by $a$,   by a vertex
disjoint collection of complete bipartite graphs. When the degrees of the vertices on one side are  $a_1,\ldots, a_b$, with average degree $a$, the following upper bound holds: $\Pi_{i} (a_i!)^{1/a_i}  \le (a!)^{b/a}$. 

{\bf Van der Waerden  inequality:} \  In 1926 Van der Waerden~\cite{waer} conjectured that the minimum
possible value of the permanent of a $b$-by-$b$  
doubly stochastic matrix is $b!/b^b$,
achieved by the matrix in which all entries are $1/b$. Proofs
of this conjecture were given in the early 80s by
Falikman~\cite{falikman} and by Egorychev~\cite{egor},
see also  Gyires~\cite{gyires}.

{\bf Schrijver's bound:}\ Schrijver~\cite{Sch98} (see also \cite{Gu} for subsequent work)
proved that every $a$-regular
bipartite graph with $b$ vertices in each class has at least
\vglue -0.3cm
$$\left( \frac{(a-1)^{a-1}}{a^{a-2}}\right)^b$$
\vglue -0.2cm
\noindent  perfect
matchings.
\vspace{0.2cm}

\noindent
{\bf Proof of Proposition~\ref{perfectmatching} (i).}
The  Bregman-Minc inequality  gives that 
$$m(G)\le
(a!)^{b/a}\le [\sqrt{2\pi a}\cdot \left(\frac{a}{e}\right)^a\cdot
e^{1/(12a)}]^{b/a} = 2^{(1/2+o(1))\log a\cdot (b/a)}\cdot \left(\frac{a}{e}\right)^b.$$
The Van der Waerden  inequality gives the following  lower bound on $m(G)$:
 $$m(G)\ge \frac{b!}{b^b}\cdot a^b\ge
  \sqrt{2\pi b}\cdot  \left(\frac{b}{e}\right)^b\cdot b^{-b} \cdot a^b.$$ 
  Using Schrijver's bound, and some delicate estimates on $e$, we obtain the following
improved asymptotics for the lower bound: 
\begin{align} \nonumber m(G)\ge & \left(
\frac{(a-1)^{a-1}}{a^{a-2}}\right)^b = \left(\frac{a}{e}\right)^b \left(
\frac{e\cdot (a-1)^{a-1}}{a^{a-1}}\right)^b \\ \nonumber \ge  & \ \left(\frac{a}{e}\right)^b\cdot
\left(\left(1+\frac{1}{a-1}+ \frac{1}{2(a-1)^2} + \frac{1}{6(a-1)^3}
\right) \left(1-\frac{1}{a} \right)\right)^{(a-1)b}\\ \nonumber =  & \
\left(\frac{a}{e}\right)^b\cdot \left(1+ \frac{1}{2a(a-1)} + \frac{1}{6a(a-1)^2}  
\right)^{(a-1)b}\ge \left(\frac{a}{e}\right)^b\cdot e^{b/2a}.  
\end{align}

The inequality in the second line above  is equivalent to the 
statement that
$$
e^{1/(a-1)} \geq 1+\frac{1}{a-1}+\frac{1}{2(a-1)^2}+\frac{1}{6(a-1)^3},
$$ 
which follows from the fact that the right-hand-side is the sum of the
four first terms in the power series of the left-hand-side, in which all
terms are positive. The final inequality is equivalent to the fact that
$$
e^{1/(2a(a-1))} \leq 1+\frac{1}{2a(a-1)}+\frac{1}{6a(a-1)^2}.
$$
This inequality holds for every $a \geq 2$ since the power series of the
left-hand-side is 
$$
1+\frac{1}{2a(a-1)}+\frac{1}{8a^2(a-1)^2} +\ldots
$$
The first two terms here are equal to the first two terms in the
expression above, and it is easy to see that for $a \geq 2$ the sum
of all the remaining terms is smaller than 
$\frac{1}{6a(a-1)^2}.$
\vspace{0.2cm}

\noindent
(ii) Follows from (i) by setting 
$a:=d, b:=n$ and $N=\left(\frac{d}{e}\right)^n.$
 \hfill\qed\\


\noindent 
{\bf Proof of Proposition~\ref{matching}.}
(i) Note that the second equation instantly follows from \eqref{prp13i}, we just need to prove the first equation.

To prove the lower bound, set  $t=b/\sqrt a$.
For every perfect matching $M$ of $G$, let $F$ be a random
subset  of $t$ edges of $M$, chosen uniformly among all 
subsets of cardinality $t$ of $M$. Note that $M-F$ is a matching
(of size $b-t$). This provides 
$m(G) \cdot {b \choose t}$ matchings, but the same matching may
be obtained multiple times. More precisely, the number of times such
a matching $M-F$ appears is exactly the number of perfect
matchings in the induced subgraph of $G$ on the set of vertices
$V(F)$ covered by the edges of $F$. The expected number of edges
in this induced subgraph is exactly
$$t+(b\cdot a-b)\frac{{t \choose 2}}{{b \choose 2}} < t+b.$$
Indeed, the subgraph contains exactly $t$ edges of $M$, and each
edge that does not belong to $M$ lies in this induced subgraph
with probability ${{t \choose 2}}/{{b \choose 2}} < t^2/b^2=1/a$. 
The above estimate thus follows from the linearity of expectation.
By Markov's Inequality it follows that  with probability  at least,
say, $1/3a$, the number of edges in this induced subgraph is smaller
than $b+2t$. This gives that with probability at least $1/3a$ the induced
subgraph of $G$ on $V(F)$ has average degree smaller than
$(b+2t)/t=\sqrt a+2$. By Minc Conjecture and the fact that
as shown in \cite{Al0}, Corollary~2.3 (see also \cite{as}, 
page 66) the upper bound provided by the Bregman-Minc 
inequality for the permanent of
a $\{0, 1\}$-matrix with a given number of $1$ entries is obtained when
all the rows have the same sum, it follows that 
the number of perfect matchings spanned by $V(F)$ is at most
$$[(\sqrt a+2)!]^{t/(\sqrt a+2)} \leq [(1+o(1))\sqrt {a}/e]^{t}.$$

Therefore, in each of these cases the matching $M-F$
obtained is counted at most 
$[(1+o(1))\sqrt {a}/e]^{t}$
times. It follows that the number of matchings of size $b-t$
in $G$ is at least
$$
\frac{1}{3a}\cdot  m(G)\cdot  \frac{{b \choose t}}{[(1+o(1))\sqrt {a}/e]^{t}}
=m(G)\cdot [(1+o(1))e^2]^{b/\sqrt a},$$
implying the desired lower bound.

For the proof of the upper bound for $m'(G)$ 
denote by $M(t)$  the set of matchings in $G$ with $b-t$ edges,
where $0\le t\le b$.
For a given $t$, there are $\binom{b}{t}^2$ ways to choose the set
of uncovered vertices, and given this choice, by the Bregman-Minc
inequality, there are at most $(a!)^{(b-t)/a}$  ways to place a perfect
matching on the rest of the vertices. Therefore 

\vglue -0.3cm
$$
|M(t)|\le \binom{b}{t}^2\cdot
(a!)^{(b-t)/a}.
$$
\vglue -0.1cm

To bound the sum of these terms for all $t$ 
define $x=(a!)^{1/a}$ and observe that this sum satisfies
\vglue -0.4cm
$$
\sum_t {b \choose t}^2 x^{b-t} \leq 
\left( \sum_t {b \choose t} (\sqrt x)^{b-t}\right)^2 =(1+\sqrt x)^{2b} = (\sqrt x)^{2b} \left( 1+\frac{1}{\sqrt x}\right)^{2b}.
$$
\vglue -0.2cm
\noindent Using Stirling's formula it is easy to check that
$$
(\sqrt x)^{2b} =(a!)^{b/a} \leq \left(\frac{a}{e}\right)^b e^{O(b \log a/a)}\quad \quad \text{and} \quad \quad 
\left(1+\frac{1}{\sqrt x}\right)^{2b} \leq e^{O(b/\sqrt x)} =e^{O(b/\sqrt a)}.
$$
This provides the required upper bound
$\left(\frac{a}{e}\right)^b \cdot 2^{O(b/\sqrt a)}.$

\noindent (ii) Follows from (i) by setting $a:=d,\  b:=n$ and 
$N=\left(\frac{d}{e}\right)^n.$
\hfill\qed\\
\vspace{0.2cm}

\noindent
Next we describe the proof of Proposition ~\ref{expgap}, starting with
an outline of the proof.
We have to prove that there is a constant $c>1$
such that $f(d)> c^n\cdot m'(d)$.  To do so we  fix a small positive
constant $\alpha$ and choose randomly and uniformly $(1+o(1))\alpha n$
$2$-dimensional subcubes of $Q_d$. Whenever two such chosen  cubes
have common vertices, we remove one of them, noting that typically the
number of subcubes removed is only $O(\alpha^2)n$. This gives a collection
$B$ of some $\beta n=(\alpha-O(\alpha^2)) n$ 
$2$-dimensional subcubes, and by lower bounding the probability that
the set $B$ produced is of size $\beta n$ we get that there are many 
distinct choices for such sets. 

We complete this partial covering by placing a
(nearly) perfect matching on the
rest of the graph, which has $(1-2\beta)n$ vertices in each of
the vertex classes, and is roughly $(1-2\beta)d$ regular. 
Since the rest of the graph is not
exactly regular, we do not have a good lower bound on the number
of its matchings. Therefore, we create a regular bipartite
graph by adding some vertices and removing some edges to make it regular.
This enables us to apply the lower bound in Proposition 
\ref{perfectmatching} and to estimate the number of ways to add a matching
and single vertices ($0$-dimensional cubes) for each collection
$B$ of $2$-cubes. A careful computation then yields the desired 
estimate. 
\vspace{0.02cm}

\noindent
\noindent {\bf Proof of Proposition~\ref{expgap}.}
The goal  is to prove that there is a constant $c>1$
such that $f(d)> c^n\cdot m'(d)$.  Let $\alpha>0$ be a sufficiently small
constant, $c>1$ will depend on the choice of $\alpha$.  Consider the
hypercube $Q_d$ as a bipartite graph with vertex classes $(U,V)$,
(the even vertices and the odd ones). 
Let $A\subset  V$ be a random subset of $V$ obtained by picking each
vertex of $V$, randomly and independently, to lie in $A$ with
probability $\alpha$. Therefore the size of $A$ is close to
$\alpha n$ with high probability. 
For each $u_i\in A$ choose a vertex $v_i\in V$ uniformly and randomly
among the $\binom{d}{2}$ choices so that $\{u_i, v_i\}$ 
together with two additional vertices from $U$ span a copy of $C_4$.

This way we get a partial covering of $Q_d$ by the $2$-cubes spanned by
$\{(u_i,v_i): \ u_i\in A\}$.
We try to complete this covering by placing a
(nearly) perfect matching on the
rest of the graph, which has about $(1-2\alpha)n$ vertices in each of
the vertex classes, and is roughly $(1-2\alpha)d$ regular. Our gain on the
number of cube partitions will come from the number of ways of choosing
the set of $2$-cubes. Unfortunately, the rest of the graph is not
exactly regular, hence we do not have a good lower bound on the number
of its matchings. To go around this we create a regular bipartite
graph by adding some vertices and removing some edges to make it regular.

In case two $2$-cubes are overlapping, we remove one of them.  
With high
proba\-bi\-lity, the number of removed vertices is $\Theta(\alpha^2 n)$. 
We thus assume this is the case, (and count 
only partitions in which this holds).
This way we obtain a partial $C_4$-covering  
$B$ of $Q_d$ (that is, $B$ is a
family of vertex disjoint $C_4$'s), with $|B|=(\alpha- \Theta(\alpha^2))
n=:\beta n,$ i.e., $(\alpha-\beta)n =\Theta(\alpha^2) n$.
When $\alpha>0$ is sufficiently small, then the number
of choices for $B$ is at least

\begin{equation}\label{beta1}
\frac{1}{n^2}\binom{n}{\alpha n} \cdot \binom{d}{2}^{\alpha n} \cdot  
\binom{n}{(\alpha-\beta )n}^{-1} \cdot \binom{d}{2}^{-(\alpha-\beta ) n } 
\ge 2^{(h(\alpha)/2)n}\cdot {d}^{2\beta n}.
\end{equation}

The first $1/n^2$ factor above is for considering only random
choices of $A$ in which $|A|=\alpha n$ and $|B|=\beta n$, meaning that
$B$ is obtained from $A$ by removing exactly $(\alpha-\beta)n$ $2$-cubes.
There are at least 
$\frac{1}{n^2}\binom{n}{\alpha n} \cdot \binom{d}{2}^{\alpha n}$
ways to choose the collection $A$ so that this holds, and each such $A$
produces a collection $B$. The number of times a fixed collection $B$
can be obtained this way 
is the number of ways to add $(\alpha-\beta)n$ $2$-cubes to this fixed 
collection, where these added subcubes contain  $(\alpha-\beta)n$ distinct
vertices $u_i$ in $V$ and the subcube is determined by $u_i$ and a vertex 
$v_i$ of Hamming distance $2$ from $u_i$. There are at most 
$\binom{n}{(\alpha-\beta )n}$ ways to choose the vertices $u_i$ and
given those, at most 
$\binom{d}{2}^{(\alpha-\beta ) n}$ to choose the corresponding vertices 
$v_i$. Therefore, dividing by the product
$\binom{n}{(\alpha-\beta )n} \cdot \binom{d}{2}^{(\alpha-\beta ) n}$
ensures that each partial collection $B$ of pairwise disjoint 
$2$-subcubes is 
counted at most once.
The last inequality follows from the fact that for small fixed $\alpha>0$
$$
{n \choose {\alpha n}}\  =\ 2^{(1+o(1)) h(\alpha)n}\gg
n^2 2^{\alpha n} \cdot 2^{h(O(\alpha^2))n}\ =\ n^2 
 \cdot 2^{\alpha n} \cdot {n \choose {(\alpha-\beta)n}} . 
$$

Let $H$ be the (random) graph spanned by $Q_d-V(B)$. Observe
that if $u,v\in V(H)$ and the distance between them in $Q_d$ 
is a least, say,
$10$, then their degrees in this graph are 
independent.


\begin{claim}\label{degree}
With high probability, all but at most 
$n/d$ vertices in $H$ have degree
in  the interval $J:= [(1- 2\beta)d-d^{2/3}, (1- 2\beta)d+d^{2/3}]$.
\end{claim}


{\bf Proof.} 
It is easy to see that for each fixed vertex 
of $H$, its degree in $H$
lies in the interval $J$ with high probability. This follows, for
example, from Azuma's martingale inequality (c.f., e.g.,
\cite{as}).
Note, however, that this does not suffice to imply the claim,
as the events corresponding to distinct vertices are not independent.  

To complete the proof of the claim partition 
$V(H)$ into $d^{10}$ classes, so that in each class the vertices
are at distance at least $10$ from each other. In each class, the
degrees of the vertices in $H$ are independent. Applying
the Chernoff bound to each class, we obtain that each class contains about
the expected number of vertices whose degrees are not in the interval.
The claimed result follows by the union bound, with room to
spare.
\hfill\qed

\smallskip
Returning to the proof of the proposition, we
assume that the assertion of the claim holds (and only count such
partitions). We thus assume  that there are at most
$n/d$ vertices in $H$ with degrees not in $J$.
As long as there is a vertex with 
degree larger than $(1- 2\beta)d+d^{2/3}$, we
remove an arbitrary subset of its edge set to make its degree $(1-
2\beta)d+d^{2/3}$.  This way we remove a total of 
at most $n/ d\cdot 2\beta\cdot
d=2\beta  n$ edges.  The number of edges missing in each vertex
class in order to get a  $((1-
2\beta)d+d^{2/3})$-regular graph is smaller than
$$n/d\cdot ((1- 2\beta)d+d^{2/3})
+ n\cdot 2  d^{2/3} + 2\beta  n < 2.5 d^{2/3}\cdot n.$$

We now add auxiliary vertices to $H$, to create a  $((1-
2\beta)d+d^{2/3})$-regular graph $F$. This can be done by 
adding at most $2.5 d^{2/3}\cdot n/(d(1-2 \beta)) \ 
\leq \  3 d^{-1/3}\cdot n$
vertices to each class.

Consider an arbitrary fixed perfect matching $M$ of $F$. 
By the lower bound in Proposition 
\ref{perfectmatching}, the number of such $M$ is at least
\vglue -0.4cm

$$\left( \frac{(1- 2\beta)d+d^{2/3}}{e}\right)^{(1-2\beta )n}\ge
(1- 2\beta)^{(1- 2\beta)n}\cdot \left(\frac{d}{e}\right)^{(1-
2\beta)n}.$$

Using $M$, we define the following cube partition of $Q_d-B$.
We let the edges of $M$ be the $1$-dimensional cubes,
with the (obvious) restriction that if an 
endpoint of an edge of $M$ is not in
$Q_d$, then it is not used, and if exactly one end point of the edge
is  in $V(Q_d)$, then it is considered as a $0$-dimensional subcube
of the partition. Note that the number of such $0$-cubes is
at most $3n/d^{1/3}$ in each class.

This way, we obtain a cube covering using subcubes of dimensions $0,1,2$;
but we might obtain the same covering from different matchings $M$. The
overcounting is at most the number of ways a perfect matching could
be placed on the vertices of $V(F)-V(Q_d)$ and the vertices of the
$0$-dimensional cubes. The number of such vertices in $V$ is at most
$6n/d^{1/3}$. Hence, by the  Bregman-Minc inequality,
the number of ways to cover them with a perfect matching is at most
$(d!)^{6n/d^{4/3}}$.  To summarize, the number of  tilings using subcubes
of dimensions $\{0,1,2\}$ is at least
\vglue -0.4cm

\begin{equation}\label{beta2}
\nonumber \frac{1}{2} \cdot 
2^{(h(\alpha)/2)n}\cdot {d}^{2\beta n}\cdot (1- 2\beta)^{(1-
2\beta)n}\cdot \left(\frac{d}{e}\right)^{(1- 2\beta)n}
\cdot(d!)^{-6n/d^{4/3}}> 2^{(h(\alpha)/3)n}\cdot N >   c^n \cdot
m'(d).
\end{equation}
The first $1/2$ factor takes care of the fact that we count only
partitions in which the required high probability events hold. 
The constant $c$ can be chosen, for example, to be
$c= 2^{(h(\alpha)/4)}>1$, where $\alpha>0$ is a 
sufficiently small absolute constant. \hfill\qed\\


\section{Upper bound on $f(d)$ and related estimates}
\label{uppfd}
\vspace{0.2cm}

\noindent {\bf Proof of Proposition~\ref{main}.}
We need to  show that
$f(d) \leq (d+1)^n$. 
Let $v_1,v_2, \ldots ,v_n$ be an enumeration of all 
the $n$ even vertices of the hypercube $Q_d$ (for example, 
lexicographically), and fix a similar enumeration of all
odd vertices. 
For each partition $P$ of the hypercube into subcubes we
construct a sequence $S=S(P)=(s_1,s_2, \ldots ,s_n)$ 
of length $n$ over the 
alphabet $\{0,1,\ldots ,d\}$
so that $S(P)$ completely defines the partition $P$.  This is done as
follows. Each element $s_i$ of $S(P)$ corresponds to the vertex $v_i$.
If $v_i$ lies in a subcube $D$ of the partition $P$ of positive dimension,
and it is a
neighbor of the lexicographically first odd vertex $u$ of
$D$, then $s_i=j$, where $j$ is the coordinate in which
$v_i$ and $u$ differ. If $v_i$ lies in such a subcube $D$ and 
is not a neighbor of $u$ then $s_i=0$, and this is also the case if
$v_i$ forms a subcube of $P$ of dimension $0$. This completes the 
definition of $S=S(P)$. It is easy to see that given $S=S(P)$, we can 
reconstruct all the odd vertices which are the lexicographically smallest
ones in subcubes in the partition $P$ of positive dimension. For each
of those, we can reconstruct the corresponding subcubes, and then we also
get all the remaining vertices which are subcubes of dimension $0$.
This shows that $f(d) \leq (d+1)^n$ and completes the proof of the 
proposition.
\hfill\qed
\smallskip

\noindent {\bf Proof of Proposition~\ref{new}.} 
 First we prove the lower bound.
A slight modification
of the Frankl-R\"odl~\cite{frankl} nibble method, as pointed out by
Grable and Phelps~\cite{grable}, gives an estimate on the number of
almost perfect matchings in $r$-uniform hypergraphs. We use a version
from Asratian and Kuzjurin~\cite{asratin}: 
\begin{theorem}
\label{asr}
The following holds for every fixed small $\delta>0$.
Let $r$ be fixed and let $\mathcal{H}_{2n}$ be an  $r$-uniform $D$-regular
hypergraph on $2n$ vertices, where both $D$ and $n$ are tending to
infinity. Furthermore, assume that the maximum codegree is $o(D)$,
where the codegree of a pair of vertices is the number of hyperedges
containing both. Then, the number of matchings
of $\mathcal{H}_{2n}$ with at least $(1-\delta)2n/r$ hyperedges is at
least $D^{(1-2\delta)(2n/r)}$.  \end{theorem}

To prove the proposition, define
a $2^r$-uniform hypergraph $\mathcal{H}_{2n}$ on the vertex
set of $Q_d$, where the hyperedges are the $r$-subcubes of $Q_d$. The
hypergraph $\mathcal{H}_{2n}$   is $D=\binom{d}{r}$-regular, with maximum
codegree $\binom{d-1}{r-1}$. A matching $M$ of size $(1-\delta)n/2^{r-1}$ corresponds to a
$\{0,r\}$-covering of $Q_d$, where if a vertex is not covered by $M$, then
we  cover it by a $0$-cube, and the hyperedges of $M$ correspond
to $r$-cubes. As $\delta$ can be chosen to be arbitrarily small,
this implies that
$$f_{0,r}(d)\ \ge\ \binom{d}{r}^{(2^{1-r}-o(1))n } = N^{r/2^{r-1}+o(1) }.$$
  
To prove the upper bound we apply the same encoding that appears
in the proof of Proposition~\ref{main} to each of the relevant partitions
here. The improved bound is obtained since every sequence that appears
in this encoding contains only a small number of nonzero
entries. Indeed, note
that each partition here contains no subcubes of dimensions
$1, 2, \ldots, r-1$, 
and for each subcube $D$
of dimension $k>0$ (which must be at least $r$) we get exactly
$k$ nonzero elements of the sequence, and $2^{k-1}-k$ zeros. All other
elements of the sequence that correspond to subcubes of dimension $0$
(consisting of an even vertex) are also $0$. 
This gives sequences of length $n$ in which the fraction of non-zeros is
at most $r/2^{r-1}$ providing an upper bound of
$$\sum_{j \leq rn/2^{r-1}} {n \choose j} d^{j}= N^{r/2^{r-1}+o(1)},$$ that provides the 
required result.
\hfill\qed

\smallskip

\noindent
{\bf Proof of Proposition~\ref{approx}}.
We estimate the number of cube partitions which do
not use  $2$-dimensional subcubes. To do so 
we bound the number of partitions in which
exactly $T$ vertices 
in each of the two vertex classes, the even ($U$) and the odd ($V$),
are covered by cubes of dimension at least $3$,
 where here $0\le T\le n$:\\
-- there are at most 
$\binom{n}{T}$ ways to choose a set $A$ of  $T$ even vertices covered 
by cubes of dimension at least $3$,\\
-- given such  choices, we can construct a sequence of length 
$T$ over the alphabet $\{0,3,4,\ldots ,d\}$ by following the procedure
in the previous proofs, but here it is applied only to the set
$A$ of $T$ chosen
even vertices. More precisely, the sequence produced corresponds to
the enume\-ration of the vertices of $A$ induced by the lexicographic 
order.  For each subcube $D$ of the partition of dimension 
at least $3$ let $u$ be its lexicographically first odd vertex.
If a vertex $v$ of $A$ belongs to $D$ and is a neighbor of $u$ then the
corresponding element of the sequence is $j$, where $u$ and $v$
differ in coordinate $j$, and if $v$ is in $D$ but is not a neighbor of
$u$ then the corresponding element is $0$. This produces 
a sequence of length $T$ over $\{0,1,2, \ldots ,d\}$
in which the fraction
of nonzero elements is at most $3/4$. Note that the sequence enables one
to reconstruct all the subcubes of dimension at least $3$ 
in the partition. Therefore, for each fixed choice of $A$ there are
less than $2^T\cdot d^{3T/4}$ ways to
place the cubes of dimension at least $3$.\\
-- Using Proposition~\ref{matching} (the statement holds for maximum degree $a$, instead of $a$-regular) with $a=d$ and $b=n-T$, given the above choices, there are at most 
$\left(\frac{d}{e}\right)^{n-T} \cdot 2^{\Theta((n-T)/\sqrt d)}$
ways to place a perfect matching on
the rest of the vertices.\\
Hence, for a  fixed $T$ the number of these cube partitions is at most
$$\binom{n}{T}  \cdot 2^T\cdot d^{3T/4}\cdot \left(\frac{d}{e}\right)^{n-T} \cdot 2^{\Theta((n-T)/\sqrt d)}.$$
  Summing up over all choices of
 $T$, and using the binomial theorem,  we obtain the upper bound 
 $$(1+2ed^{-1/4})^n \cdot 2^{\Theta(n/\sqrt d)} \cdot N \le \exp(20 n/d^{1/4})\cdot N.$$ \hfill\qed

\noindent{\bf Proof of Proposition~\ref{ired}.}
The construction is recursive. By induction on the dimension $d$
we construct many irreducible tight partitions of $Q_d$, all containing
one specific subcube. We start in some fixed dimension $d$ with
at least $3$ irreducible tight partitions of $Q_d$ all containing
some subcube. Given
two distinct partitions $B_1,B_2$ among those 
consider the partition $B =B_10 \cup B_21$ of $Q_{d+1}$ in which
the partition $B_1$ is used in the first hyperplane in which the last
coordinate is $0$ and the partition $B_2$ is used in the second
hyperplane in which the last coordinate is $1$. Since
$B_1$ and $B_2$ have the same subcube $D$ we replace 
$D0$ and $D1$ in $B$ by their union $D'=D_0 \cup D_1$.
By replacing all (the nonzero number) of such $D$
we obtain this way an irreducible tight partition of $Q_{d+1}$ from each
ordered pair of partitions $B_1,B_2$ we had, and all these new partitions
contain the same subcube $D'$. It is not difficult to see that all the
partitions obtained are tight since
$B_1,B_2$ are tight and different. Each partition obtained is also
irreducible since no subcube in which the last coordinate
is fixed can be spanned by a subfamily, as $B_1,B_2$ are
irreducible. Similarly, if 
a subcube in which the last coordinate is not fixed is
spanned by a subfamily, then the construction ensures that its
two subcubes consisting of all vertices with a fixed last coordinate 
are also spanned by a subfamily of $B_1$ and a subfamily of $B_2$.
As all the partitions obtained are distinct,
this shows that if the number of
irreducible tight partitions we get for dimension $d$ in this way is
$x_d$, then $x_{d+1} =x_d (x_d-1)$, implying the desired lower bound.
\vspace{0.2cm}

\noindent
{\bf Acknowledgment.}
We are grateful to Maximilien Gadouleau for pointing out
reference~\cite{BGT20}, and for 
explaining the motivation.  We also thank Artur
Riazanov and Yuval Filmus for useful comments. Part of this work was
completed while the second author visited Tel Aviv University, 
partially supported by the Simons Foundation,  by the Mathematisches
Forschungsinstitut Oberwolfach, and by the ERC Consolidator Grant
101044123 (RandomHypGra) of Wojciech Samotij.
Finally, we thank the anonymous referees for many useful comments.

\end{document}